\documentclass [12pt]{article}
\usepackage{graphicx,amssymb,amsfonts,latexsym,amsmath,amsthm,times}
\usepackage{epsfig}
\usepackage{color}
\usepackage[english]{babel}
\usepackage[a4paper]{geometry}
\def\lanbox{\hbox{$\, \vrule height 0.25cm width 0.25cm depth 0.01cm \,$}}

\numberwithin{equation}{section}

\begin{document}

\vspace*{1.4cm}

\normalsize \centerline{\Large \bf ON THE WELL-POSEDNESS
OF A CERTAIN MODEL WITH THE }

\medskip

\centerline{\Large \bf BI-LAPLACIAN APPEARING IN THE MATHEMATICAL BIOLOGY}

\vspace*{1cm}

\centerline{\bf Messoud Efendiev$^{1,2}$, Vitali Vougalter$^{3 \ *}$}

\bigskip

\centerline{$^1$ Helmholtz Zentrum M\"unchen, Institut f\"ur Computational
Biology, Ingolst\"adter Landstrasse 1}

\centerline{Neuherberg, 85764, Germany}

\centerline{e-mail: messoud.efendiyev@helmholtz-muenchen.de}

\centerline{$^2$ Department of Mathematics, Marmara University, Istanbul,
T\"urkiye}

\centerline{e-mail: m.efendiyev@marmara.edu.tr}

\bigskip

\centerline{$^{3 \ *}$  Department of Mathematics, University
of Toronto}

\centerline{Toronto, Ontario, M5S 2E4, Canada}

\centerline{ e-mail: vitali@math.toronto.edu}

\medskip


\vspace*{0.25cm}

\noindent {\bf Abstract:}
The work is devoted to the global well-posedness in
$W^{1,(4,2)}({\mathbb R}\times {\mathbb R}^{+})$ of the
integro-differential problem involving the square of the 
one dimensional Laplace operator along with the drift
term. Our proof is based on a fixed point technique. 
Moreover, we provide the assumption leading to the
existence of the nontrivial solution for the problem
under the consideration.
Such  equation is relevant to the cell
population dynamics in the Mathematical Biology.

\vspace*{0.25cm}



\vspace*{0.5cm}

\bigskip

\bigskip


\setcounter{section}{1}

\centerline{\bf 1. Introduction}

\medskip

\noindent
The present article is devoted to the global well-posedness
of the nonlocal reaction-diffusion equation with the constants
$\displaystyle{a\geq 0}$ and $b\in {\mathbb R}$, namely
\begin{equation}
\label{h}
\frac{\partial u}{\partial t} =
-\frac{\partial^{4}u}{\partial x^{4}}+
b\frac{\partial u}{\partial x}+au+
\int_{-\infty}^{\infty}G(x-y)F(u(y,t), y)dy, \quad x\in {\mathbb R}
\end{equation}
relevant to the cell population dynamics. We assume that the initial
condition for (\ref{h}) is 
\begin{equation}
\label{ic}
u(x,0)=u_{0}(x)\in H^{4}({\mathbb R}).
\end{equation} 
The similar problem on the real line involving the  fractional Laplacian in the context of the anomalous diffusion was treated
in ~\cite{EV25}.
Note that the existence of stationary solutions of the integro-differential equations 
with the bi-Laplacian and the biological applications of such models but without a transport term were discussed in ~\cite{VV210}.
The cases on the whole real line and
on a finite interval with periodic boundary conditions involving the drift term and the square root of the one dimensional
negative Laplacian were covered in
~\cite{EV22}. The work  ~\cite{EV20} deals with situation of the normal diffusion and the transport.
Solvability of certain integro-differential equations with anomalous diffusion, transport
and the cell influx/efflux was considered in ~\cite{VV21}.
Spatial structures and generalized travelling waves for an
integro-differential equation were treated in ~\cite{ABVV10}. Spatial patterns appearing in higher order models in physics
and mechanics were covered in ~\cite{PT01}. The article
~\cite{VV130} is devoted to the 
emergence and propagation of patterns in nonlocal reaction-
diffusion equations arising in the theory of speciation and containing the
drift term.  Pattern and waves for a model in population dynamics 
with nonlocal consumptions of resources were studied in ~\cite{GVA06}.
The existence of steady states and travelling waves for the
non-local Fisher-KPP equation was covered in ~\cite{BNPR09}. In
~\cite{BHN05} the authors estimated the speed of propagation for KPP type
problems in the periodic framework. Important applications to the theory of
reaction-diffusion equations with non-Fredholm operators were developed in
~\cite{DMV05}, ~\cite{DMV08}. Fredholm structures, topological invariants and applications were covered in ~\cite{E09}.
Evolution equations arising in the modelling of life sciences were considered in ~\cite{E13}.
The work ~\cite{GK18} deals with the entropy method for generalized Poisson-Nernst-
Planck equations. The large time behavior of solutions of fourth order parabolic equations and $\epsilon$-entropy of their attractors were discussed in ~\cite{EP07}. Lower estimate of the attractor dimension for a chemotaxis growth system was performed in
~\cite{ATEYM06}. The theory of finite and infinite dimensional attractors for evolution equations of mathematical physics was
developed in ~\cite{E10}. Attractors for degenerate parabolic type equations were investigated in ~\cite{E131}.
Exponential decay toward equilibrium via entropy methods 
for reaction-diffusion equations was established in ~\cite{DF06}.  Local and global existence for nonlocal multispecies 
advection-diffusion models were covered in ~\cite{GHLP22}.
A free boundary problem driven by the
biharmonic operator was studied in ~\cite{DKV20}.

The space variable $x$ in our article is correspondent to the cell genotype,
$u(x,t)$ denotes the cell density as a function of the genotype and time.
The right side of (\ref{h}) describes the evolution of the cell density via
the cell proliferation, mutations and transport.
The diffusion term corresponds to the change of genotype
due to the small random mutations, and the integral term describes large mutations.
The function $F(u, x)$ stands for the rate of cell birth depending on $u$ and $x$
(density dependent proliferation), and the kernel $G(x-y)$ gives
the proportion of newly born cells which change their genotype from $y$ to $x$.
Let us assume that it depends on the distance between the genotypes.

The standard Fourier transform in our context equals to
\begin{equation}
\label{ft}
\widehat{\phi}(p):=\frac{1}{\sqrt{2\pi}}\int_{-\infty}^{\infty}\phi(x)e^{-ipx}dx,
\quad p\in {\mathbb R}.
\end{equation}
Clearly, the estimate from above
\begin{equation}
\label{fub}  
\|\widehat{\phi}(p)\|_{L^{\infty}({\mathbb R})}\leq \frac{1}{\sqrt{2\pi}}
\|\phi(x)\|_{L^{1}({\mathbb R})}
\end{equation}
is valid (see e.g. ~\cite{LL97}). Evidently, (\ref{fub}) implies that
\begin{equation}
\label{fub1}  
\|p^{4}\widehat{\phi}(p)\|_{L^{\infty}({\mathbb R})}\leq \frac{1}{\sqrt{2\pi}}
\Big\|\frac{d^{4}\phi}{dx^{4}}\Big\|_{L^{1}({\mathbb R})}.
\end{equation}
Let us suppose that the conditions below on the integral kernel contained in 
equation (\ref{h}) are satisfied.

\medskip

\noindent
{\bf Assumption 1.1.}  {\it Let $G(x): {\mathbb R}\to {\mathbb R}$ be
nontrivial, such that
$\displaystyle{G(x), \ \frac{d^{4}G(x)}{dx^{4}}\in L^{1}({\mathbb R})}$.}     

\medskip

\noindent
This allows us to define the auxiliary expression
\begin{equation}
\label{g}  
g:=\sqrt{\|G(x)\|_{L^{1}({\mathbb R})}^{2}+
\Big\|\frac{d^{4}G(x)}{dx^{4}}\Big\|_{L^{1}({\mathbb R})}^{2}}.
\end{equation}
Thus, $0<g<\infty$.

\medskip

\noindent
From the perspective of the applications, the space dimension is
not limited to $d=1$ since the space variable corresponds to the cell
genotype but not to the usual physical space.
We have the Sobolev space
\begin{equation}
\label{ss}  
H^{4}({\mathbb R}):=
\Bigg\{\phi(x):{\mathbb R}\to {\mathbb {\mathbb R}} \ | \
\phi(x)\in L^{2}({\mathbb R}), \ \frac{d^{4}\phi}{dx^{4}}\in
L^{2}({\mathbb R}) \Bigg \}.
\end{equation}
It is equipped with the norm
\begin{equation}
\label{n}
\|\phi\|_{H^{4}({\mathbb R})}^{2}:=\|\phi\|_{L^{2}({\mathbb R})}^{2}+
\Bigg\| \frac{d^{4}\phi}{dx^{4}}\Bigg\|_{L^{2}({\mathbb R})}^{2}.
\end{equation}
For establishing the global well-posedness of problem (\ref{h}), (\ref{ic}),
we will use the function space
$$
W^{1, (4, 2)}({\mathbb R}\times [0, T]):=
$$
\begin{equation}
\label{142}
\Big\{u(x,t): {\mathbb R}\times [0, T]\to {\mathbb R} \ \Big|
\ u(x,t), \ \frac{\partial^{4}u}{\partial x^{4}}, \
\frac{\partial u}{\partial t}
\in L^{2}({\mathbb R}\times [0, T]) \Big\},
\end{equation}  
so that
$$
\|u(x,t)\|_{W^{1, (4, 2)}({\mathbb R}\times [0, T])}^{2}:=
$$
\begin{equation}
\label{142n}
\Big\|\frac{\partial u}{\partial t}\Big\|_{L^{2}({\mathbb R}\times [0, T])}^{2}+
\Big\|\frac{\partial^{4}u}{\partial x^{4}}\Big\|_{L^{2}({\mathbb R}\times [0, T])}^{2}+
\|u\|_{L^{2}({\mathbb R}\times [0, T])}^{2},
\end{equation}  
where $T>0$. In definition (\ref{142n}) we have
$$
\|u\|_{L^{2}({\mathbb R}\times [0, T])}^{2}:=\int_{0}^{T}\int_{-\infty}^{\infty}
|u(x,t)|^{2}dxdt.
$$
Throughout the work we will also use the norm
$$
\|u(x,t)\|_{L^{2}({\mathbb R})}^{2}:=\int_{-\infty}^{\infty}|u(x,t)|^{2}dx.
$$

\medskip

\noindent
{\bf Assumption 1.2.} {\it Function
$F(u, x): {\mathbb R}\times{\mathbb R}\to {\mathbb R}$ is satisfying the
Caratheodory condition (see ~\cite{K64}), so that
\begin{equation}
\label{lub}
|F(u, x)|\leq k|u|+h(x) \quad for \quad u\in {\mathbb R}, \quad x\in {\mathbb R}
\end{equation}  
with a constant $k>0$ and
$h(x): {\mathbb R}\to {\mathbb R}^{+}, \ h(x)\in L^{2}({\mathbb R})$.
Moreover, it is a Lipschitz continuous function, such that
\begin{equation}
\label{lip}
|F(u_{1}, x)-F(u_{2}, x)|\leq l|u_{1}-u_{2}| \quad for \quad any \quad
u_{1, 2}\in {\mathbb R}, \quad x\in {\mathbb R}
\end{equation}
with a constant $l>0$.}

\medskip

\noindent
In our article ${\mathbb R}^{+}$ stands for the nonnegative semi-axis.
The solvability of a local elliptic equation in a bounded domain in
${\mathbb R}^{N}$ was covered in ~\cite{BO86}. The nonlinear function involved there
was allowed to have a sublinear growth.

We apply the standard Fourier transform (\ref{ft}) to both sides of
problem (\ref{h}), (\ref{ic}). This gives us
\begin{equation}
\label{hf}
\frac{\partial \widehat{u}}{\partial t}=[-p^{4}+ibp+a]\widehat{u}+
\sqrt{2\pi}\widehat{G}(p)\widehat{f}_{u}(p,t),
\end{equation}
\begin{equation}
\label{icf}
\widehat{u(x,0)}(p)=\widehat{u_{0}}(p).  
\end{equation}
In formula (\ref{hf}) and below $\widehat{f}_{u}(p,t)$ will denote the Fourier image
of $F(u(x,t), x)$. Obviously,
$$
u(x,t)=\frac{1}{\sqrt{2\pi}}\int_{-\infty}^{\infty}\widehat{u}(p,t)e^{ipx}dp,
\quad 
\frac{\partial u}{\partial t}=\frac{1}{\sqrt{2\pi}}\int_{-\infty}^{\infty}
\frac{\partial \widehat{u}(p,t)}{\partial t}e^{ipx}dp,
$$
where $x\in {\mathbb R}, \ t\geq 0$. By virtue of the Duhamel's principle,
we can reformulate problem (\ref{hf}), (\ref{icf}) as
$$
\widehat{u}(p,t)=
$$
\begin{equation}
\label{duh}
e^{t\{-p^{4}+ibp+a\}}\widehat{u_{0}}(p)+\int_{0}^{t}
e^{(t-s)\{-p^{4}+ibp+a\}}\sqrt{2\pi}\widehat{G}(p)\widehat{f}_{u}(p,s)ds.
\end{equation}
Related to (\ref{duh}), we have the auxiliary equation
$$
\widehat{u}(p,t)=
$$
\begin{equation}
\label{aux}
e^{t\{-p^{4}+ibp+a\}}\widehat{u_{0}}(p)+\int_{0}^{t}
e^{(t-s)\{-p^{4}+ibp+a\}}\sqrt{2\pi}\widehat{G}(p)\widehat{f}_{v}(p,s)ds.
\end{equation}
Here $\widehat{f}_{v}(p,s)$ stands for the Fourier image of
$F(v(x,s), x)$ under transform (\ref{ft}), where $v(x,t)$ is an arbitrary function
belonging to $W^{1, (4, 2)}({\mathbb R}\times [0, T])$.

We introduce the operator $\tau_{a, b}$, so that $u =\tau_{a, b}v$, where $u$ 
solves (\ref{aux}). The main result of the article is as follows.

\bigskip

\noindent
{\bf Theorem 1.3.} {\it Let Assumptions 1.1 and 1.2 be valid and
\begin{equation}
\label{qlt}
gl\sqrt{T^{2}e^{2aT}(1+2[a+|b|+1]^{2})+2}<1 
\end{equation}
with the constant $g$ defined in (\ref{g}) and the Lipschitz constant
$l$ introduced in (\ref{lip}).
Then equation (\ref{aux}) defines the map
$\tau_{a, b}: W^{1, (4, 2)}({\mathbb R}\times [0, T])\to W^{1, (4, 2)}
({\mathbb R}\times [0, T])$, which is a strict contraction.
The unique fixed point $w(x,t)$ of this map $\tau_{a, b}$ is the only solution of
problem (\ref{h}), (\ref{ic}) in
$W^{1, (4, 2)}({\mathbb R}\times [0, T])$.}

\medskip

The final proposition of our work is devoted to the global well-posedness for
our equation. 

\medskip

\noindent
{\bf Corollary 1.4.} {\it Let the assumptions of Theorem 1.3 hold.
Then problem (\ref{h}), (\ref{ic}) possesses a unique solution   
$w(x,t)\in W^{1, (4, 2)}({\mathbb R}\times {\mathbb R}^{+})$. Such solution is
nontrivial for $x\in {\mathbb R}$ and $t\in {\mathbb R}^{+}$ 
provided the intersection of supports of the Fourier images of functions
$\hbox{supp}\widehat{F(0, x)}\cap \hbox{supp}\widehat{G}$ is a set of
nonzero Lebesgue measure on the real line.}

\medskip

Let us proceed to the proof of our main result.

\bigskip


\setcounter{section}{2}
\setcounter{equation}{0}

\centerline{\bf 2. The well-posedness of the equation}

\bigskip

\noindent
{\it Proof of Theorem 1.3.} We choose arbitrarily
$v(x,t)\in W^{1, (4, 2)}({\mathbb R}\times [0, T])$.
It can be trivially checked that the first term in the right side of (\ref{aux})
is contained in $L^{2}({\mathbb R}\times [0, T])$. Evidently,
$$
\|e^{t\{-p^{4}+ibp+a\}}\widehat{u_{0}}(p)\|_{L^{2}({\mathbb R})}^{2}=
\int_{-\infty}^{\infty}e^{-2tp^{4}}e^{2at}|\widehat{u_{0}}(p)|^{2}dp\leq
e^{2at}\|u_{0}\|_{L^{2}({\mathbb R})}^{2},
$$
so that
$$
\|e^{t\{-p^{4}+ibp+a\}}\widehat{u_{0}}(p)\|_{L^{2}({\mathbb R}\times [0, T])}^{2}=
\int_{0}^{T}\|e^{t\{-p^{4}+ibp+a\}}\widehat{u_{0}}(p)\|_{L^{2}({\mathbb R})}^{2}dt
\leq
$$
$$
\int_{0}^{T}e^{2at}\|u_{0}\|_{L^{2}({\mathbb R})}^{2}dt.
$$
The right hand side of the last inequality is equal to
$\displaystyle{\frac{e^{2aT}-1}{2a}\|u_{0}\|_{L^{2}({\mathbb R})}^{2}}$
if $a>0$ and
$T\|u_{0}\|_{L^{2}({\mathbb R})}^{2}$ for $a=0$. Thus,
\begin{equation}
\label{u0l2}  
e^{t\{-p^{4}+ibp+a\}}\widehat{u_{0}}(p)\in L^{2}({\mathbb R}\times [0, T]).
\end{equation}
Clearly, the upper bound on the norm of the second term in the right side of (\ref{aux}) 
$$
\Big\|\int_{0}^{t}
e^{(t-s)\{-p^{4}+ibp+a\}}\sqrt{2\pi}\widehat{G}(p)\widehat{f}_{v}(p,s)ds\Big\|_
{L^{2}({\mathbb R})}\leq
$$
$$  
\int_{0}^{t}\Big\|
e^{(t-s)\{-p^{4}+ibp+a\}}\sqrt{2\pi}\widehat{G}(p)\widehat{f}_{v}(p,s)\Big\|_
{L^{2}({\mathbb R})}ds.
$$
holds. We have
$$  
\Big\|
e^{(t-s)\{-p^{4}+ibp+a\}}\sqrt{2\pi}\widehat{G}(p)\widehat{f}_{v}(p,s)\Big\|_
{L^{2}({\mathbb R})}^{2}=
$$
\begin{equation}
\label{intgfv}  
\int_{-\infty}^{\infty}e^{-2(t-s)p^{4}}e^{2a(t-s)}2\pi|\widehat{G}(p)|^{2}
|\widehat{f}_{v}(p,s)|^{2}dp.
\end{equation}
Let us use inequality (\ref{fub}) to derive the estimate from above on the right side of
(\ref{intgfv}) as
$$
e^{2a(t-s)}2\pi\|\widehat{G}(p)\|_{L^{\infty}({\mathbb R})}^{2}
\|F(v(x,s),x)\|_{L^{2}({\mathbb R})}^{2}\leq
$$
$$
e^{2aT}\|G(x)\|_{L^{1}({\mathbb R})}^{2}
\|F(v(x,s),x)\|_{L^{2}({\mathbb R})}^{2}.
$$
Hence,
$$  
\Big\|
e^{(t-s)\{-p^{4}+ibp+a\}}\sqrt{2\pi}\widehat{G}(p)\widehat{f}_{v}(p,s)\Big\|_
{L^{2}({\mathbb R})}\leq
$$
$$
e^{aT}\|G(x)\|_{L^{1}({\mathbb R})}
\|F(v(x,s),x)\|_{L^{2}({\mathbb R})}.
$$
By virtue of (\ref{lub}), we obtain
\begin{equation}
\label{lubs}  
\|F(v(x,s),x)\|_{L^{2}({\mathbb R})}\leq k\|v(x,s)\|_{L^{2}({\mathbb R})}+\|h(x)\|_
{L^{2}({\mathbb R})}.
\end{equation}
Then
$$  
\Big\|
e^{(t-s)\{-p^{4}+ibp+a\}}\sqrt{2\pi}\widehat{G}(p)\widehat{f}_{v}(p,s)\Big\|_
{L^{2}({\mathbb R})}\leq
$$
$$
e^{aT}\|G(x)\|_{L^{1}({\mathbb R})}
\{k\|v(x,s)\|_{L^{2}({\mathbb R})}+\|h(x)\|_{L^{2}({\mathbb R})}\},
$$
so that
$$  
\Big\|\int_{0}^{t}
e^{(t-s)\{-p^{4}+ibp+a\}}\sqrt{2\pi}\widehat{G}(p)\widehat{f}_{v}(p,s)ds\Big\|_
{L^{2}({\mathbb R})}\leq
$$
$$
ke^{aT}\|G(x)\|_{L^{1}({\mathbb R})}\int_{0}^{T}\|v(x,s)\|_{L^{2}({\mathbb R})}ds+
Te^{aT}\|G(x)\|_{L^{1}({\mathbb R})}\|h(x)\|_{L^{2}({\mathbb R})}.
$$
By means of the Schwarz inequality
\begin{equation}
\label{sch}  
\int_{0}^{T}\|v(x,s)\|_{L^{2}({\mathbb R})}ds\leq
\sqrt{\int_{0}^{T}\|v(x,s)\|_{L^{2}({\mathbb R})}^{2}ds}\sqrt{T}.
\end{equation}
This yields
$$
\Big\|\int_{0}^{t}
e^{(t-s)\{-p^{4}+ibp+a\}}\sqrt{2\pi}\widehat{G}(p)\widehat{f}_{v}(p,s)ds\Big\|_
{L^{2}({\mathbb R})}^{2}\leq
$$
$$
e^{2aT}\|G(x)\|_{L^{1}({\mathbb R})}^{2}
\{k\sqrt{T}\|v(x,s)\|_{L^{2}({\mathbb R}\times [0,T])}+T\|h(x)\|_{L^{2}({\mathbb R})}\}^{2}.
$$
We derive the upper bound on the norm as
$$
\Big\|\int_{0}^{t}
e^{(t-s)\{-p^{4}+ibp+a\}}\sqrt{2\pi}\widehat{G}(p)\widehat{f}_{v}(p,s)ds\Big\|_
{L^{2}({\mathbb R}\times [0,T])}^{2}=
$$
$$
\int_{0}^{T}\Big\|\int_{0}^{t}
e^{(t-s)\{-p^{4}+ibp+a\}}\sqrt{2\pi}\widehat{G}(p)\widehat{f}_{v}(p,s)ds\Big\|_
{L^{2}({\mathbb R})}^{2}dt\leq 
$$
$$
e^{2aT}\|G(x)\|_{L^{1}({\mathbb R})}^{2}
\{k\|v(x,s)\|_{L^{2}({\mathbb R}\times [0,T])}+\sqrt{T}\|h(x)\|_{L^{2}({\mathbb R})}\}^{2}
T^{2}<\infty
$$
under the stated assumptions for
$v(x,s)\in W^{1, (4, 2)}({\mathbb R}\times [0, T])$. Therefore,
\begin{equation}
\label{0tetsgvl2}  
\int_{0}^{t}
e^{(t-s)\{-p^{4}+ibp+a\}}\sqrt{2\pi}\widehat{G}(p)\widehat{f}_{v}(p,s)ds\in
L^{2}({\mathbb R}\times [0,T]).
\end{equation}
By virtue of (\ref{u0l2}), (\ref{0tetsgvl2}) and (\ref{aux}), we have
\begin{equation}
\label{upt12}  
\widehat{u}(p,t)\in L^{2}({\mathbb R}\times [0,T]).
\end{equation}
This means that
\begin{equation}
\label{uxtl2}  
u(x,t)\in L^{2}({\mathbb R}\times [0,T]).
\end{equation}
Let us recall formula  (\ref{aux}). Thus,
$$
p^{4}\widehat{u}(p,t)=
$$
\begin{equation}
\label{aux2}
e^{t\{-p^{4}+ibp+a\}}p^{4}\widehat{u_{0}}(p)+\int_{0}^{t}
e^{(t-s)\{-p^{4}+ibp+a\}}\sqrt{2\pi}p^{4}\widehat{G}(p)\widehat{f}_{v}(p,s)ds.
\end{equation}
Consider the first term in the right side of (\ref{aux2}). Clearly,
$$
\|e^{t\{-p^{4}+ibp+a\}}p^{4}\widehat{u_{0}}(p)\|_{L^{2}({\mathbb R}\times [0,T])}^{2}=
\int_{0}^{T}\int_{-\infty}^{\infty}e^{-2tp^{4}}e^{2at}|p^{4}\widehat{u_{0}}(p)|^{2}dp
dt\leq
$$
$$
\int_{0}^{T}\int_{-\infty}^{\infty}e^{2at}|p^{4}\widehat{u_{0}}(p)|^{2}dpdt.
$$
Obviously, this equals to
$\displaystyle
{\frac{e^{2aT}-1}{2a}\Big\|\frac{d^{4}u_{0}}{dx^{4}}\Big\|_{L^{2}({\mathbb R})}^{2}}$ for
$a>0$ and
$\displaystyle{T\Big\|\frac{d^{4}u_{0}}{dx^{4}}\Big\|_{L^{2}({\mathbb R})}^{2}}$
when $a=0$. Therefore,
\begin{equation}
\label{p2u0hpl2}
e^{t\{-p^{4}+ibp+a\}}p^{4}\widehat{u_{0}}(p)\in L^{2}({\mathbb R}\times [0,T]).
\end{equation}  
We proceed to analyzing the second term in the right side of
(\ref{aux2}). Evidently,
$$
\Big\|\int_{0}^{t}
e^{(t-s)\{-p^{4}+ibp+a\}}\sqrt{2\pi}p^{4}\widehat{G}(p)\widehat{f}_{v}(p,s)ds
\Big\|_{L^{2}({\mathbb R})}\leq 
$$
$$
\int_{0}^{t}\Big\|
e^{(t-s)\{-p^{4}+ibp+a\}}\sqrt{2\pi}p^{4}\widehat{G}(p)\widehat{f}_{v}(p,s)
\Big\|_{L^{2}({\mathbb R})}ds.
$$
Note that
$$
\Big\|e^{(t-s)\{-p^{4}+ibp+a\}}\sqrt{2\pi}p^{4}\widehat{G}(p)\widehat{f}_{v}(p,s)
\Big\|_{L^{2}({\mathbb R})}^{2}=
$$
\begin{equation}
\label{etsp2gpl2}
\int_{-\infty}^{\infty}
e^{-2(t-s)p^{4}}e^{2a(t-s)}2\pi|p^{4}\widehat{G}(p)|^{2}|\widehat{f}_{v}(p,s)|^{2}
dp.
\end{equation}  
The right side of (\ref{etsp2gpl2}) can be estimated from above via inequality (\ref{fub1})
as
$$
2\pi e^{2aT}\|p^{4}\widehat{G}(p)\|_{L^{\infty}({\mathbb R})}^{2}
\int_{-\infty}^{\infty}|\widehat{f}_{v}(p,s)|^{2}dp\leq
$$
$$
e^{2aT}\Big\|\frac{d^{4}G}{dx^{4}}\Big\|_{L^{1}({\mathbb R})}^{2}
\|F(v(x,s),x)\|_{L^{2}({\mathbb R})}^{2}.
$$
By means of (\ref{lubs}) , we have
$$
\Big\|e^{(t-s)\{-p^{4}+ibp+a\}}\sqrt{2\pi}p^{4}\widehat{G}(p)\widehat{f}_{v}(p,s)
\Big\|_{L^{2}({\mathbb R})}\leq
$$
$$
e^{aT}\Big\|\frac{d^{4}G}{dx^{4}}\Big\|_{L^{1}({\mathbb R})}\{k\|v(x,s)\|_
{L^{2}({\mathbb R})}+\|h(x)\|_{L^{2}({\mathbb R})}\}.
$$
Then
$$
\Big\|\int_{0}^{t}e^{(t-s)\{-p^{4}+ibp+a\}}\sqrt{2\pi}p^{4}\widehat{G}(p)
\widehat{f}_{v}(p,s)ds\Big\|_{L^{2}({\mathbb R})}\leq
$$
$$
ke^{aT}\Big\|\frac{d^{4}G}{dx^{4}}\Big\|_{L^{1}({\mathbb R})}\int_{0}^{T}\|v(x,s)\|_
{L^{2}({\mathbb R})}ds+
Te^{aT}\Big\|\frac{d^{4}G}{dx^{4}}\Big\|_{L^{1}({\mathbb R})}\|h(x)\|_{L^{2}({\mathbb R})}.
$$
Let us recall estimate (\ref{sch}), which implies
$$
\Big\|\int_{0}^{t}e^{(t-s)\{-p^{4}+ibp+a\}}\sqrt{2\pi}p^{4}\widehat{G}(p)
\widehat{f}_{v}(p,s)ds\Big\|_{L^{2}({\mathbb R})}^{2}\leq
$$
$$
e^{2aT}\Big\|\frac{d^{4}G}{dx^{4}}\Big\|_{L^{1}({\mathbb R})}^{2}
\{k\|v(x,s)\|_{L^{2}({\mathbb R}\times [0,T])}\sqrt{T}+
\|h(x)\|_{L^{2}({\mathbb R})}T\}^{2}.
$$
This means that
$$
\Big\|\int_{0}^{t}e^{(t-s)\{-p^{4}+ibp+a\}}\sqrt{2\pi}p^{4}\widehat{G}(p)
\widehat{f}_{v}(p,s)ds\Big\|_{L^{2}({\mathbb R}\times [0, T])}^{2}\leq
$$
$$
e^{2aT}\Big\|\frac{d^{4}G}{dx^{4}}\Big\|_{L^{1}({\mathbb R})}^{2}
\{k\|v(x,s)\|_{L^{2}({\mathbb R}\times [0,T])}+
\|h(x)\|_{L^{2}({\mathbb R})}\sqrt{T}\}^{2}T^{2}<\infty
$$
under the given conditions with
$v(x,s)\in W^{1, (4, 2)}({\mathbb R}\times [0, T])$. Therefore,
\begin{equation}
\label{int0tp2Gpfv}  
\int_{0}^{t}e^{(t-s)\{-p^{4}+ibp+a\}}\sqrt{2\pi}p^{4}\widehat{G}(p)
\widehat{f}_{v}(p,s)ds\in L^{2}({\mathbb R}\times [0, T]).
\end{equation}
(\ref{p2u0hpl2}), (\ref{int0tp2Gpfv}) and (\ref{aux2}) yield that
\begin{equation}
\label{up2t12}  
p^{4}\widehat{u}(p,t)\in L^{2}({\mathbb R}\times [0, T]),
\end{equation}  
so that
\begin{equation}
\label{d2udx2l2}
\frac{\partial^{4}u}{\partial x^{4}}\in L^{2}({\mathbb R}\times [0, T]).
\end{equation}  
Using (\ref{aux}), we easily derive
\begin{equation}
\label{duhdt}
\frac{\partial \widehat{u}}{\partial t}=\{-p^{4}+ibp+a\}\widehat{u}(p,t)+
\sqrt{2\pi}\widehat{G}(p)\widehat{f}_{v}(p,t).
\end{equation}
According to (\ref{upt12}),
\begin{equation}
\label{aupt12}  
a\widehat{u}(p,t)\in L^{2}({\mathbb R}\times [0,T]).
\end{equation}
We obtain the upper bound for the norm as
$$
\|ibp\widehat{u}(p,t)\|_{L^{2}({\mathbb R}\times [0,T])}^{2}=b^{2}\int_{0}^{T}
\Big\{\int_{|p|\leq 1}p^{2}|\widehat{u}(p,t)|^{2}dp+
\int_{|p|>1}p^{2}|\widehat{u}(p,t)|^{2}dp\Big\}dt\leq
$$
$$
b^{2}\{\|\widehat{u}(p,t)\|_{L^{2}({\mathbb R}\times [0,T])}^{2}+
\|p^{4}\widehat{u}(p,t)\|_{L^{2}({\mathbb R}\times [0,T])}^{2}\}<\infty
$$
due to (\ref{upt12}) and (\ref{up2t12}). Hence,
\begin{equation}
\label{bupt12} 
ibp\widehat{u}(p,t)\in L^{2}({\mathbb R}\times [0,T]).
\end{equation}
And
\begin{equation}
\label{uptalf12} 
-p^{4}\widehat{u}(p,t)\in L^{2}({\mathbb R}\times [0,T])
\end{equation}
via  (\ref{up2t12}).
Let us combine statements  (\ref{aupt12}), (\ref{bupt12}) and (\ref{uptalf12}), which yields
\begin{equation}
\label{uptalfab12}
(-p^{4}+ibp+a)\widehat{u}(p,t)\in L^{2}({\mathbb R}\times [0,T]).
\end{equation}
We consider the remaining term in the right side of (\ref{duhdt}). Recall inequalities (\ref{fub}) and 
(\ref{lubs}). Clearly,
$$
\|\sqrt{2\pi}\widehat{G}(p)\widehat{f}_{v}(p,t)\|_
{L^{2}({\mathbb R}\times [0,T])}^{2}\leq 2\pi
\|\widehat{G}(p)\|_{L^{\infty}({\mathbb R})}^{2}\int_{0}^{T}
\|F(v(x,t), x)\|_{L^{2}({\mathbb R})}^{2}dt\leq
$$
$$
\|G(x)\|_{L^{1}({\mathbb R})}^{2}\int_{0}^{T}(k\|v(x,t)\|_{L^{2}({\mathbb R})}+
\|h(x)\|_{L^{2}({\mathbb R})})^{2}dt\leq
$$
$$
\|G(x)\|_{L^{1}({\mathbb R})}^{2}\{2k^{2}\|v(x,t)\|_{L^{2}({\mathbb R}\times [0, T])}^{2}+
2\|h(x)\|_{L^{2}({\mathbb R})}^{2}T\}<\infty
$$
under the stated assumptions for
$v(x,t)\in W^{1, (4, 2)}({\mathbb R}\times [0, T])$. Hence,
\begin{equation}
\label{ghfvhpt}  
\sqrt{2\pi}\widehat{G}(p)\widehat{f}_{v}(p,t)\in L^{2}({\mathbb R}\times [0,T]).
\end{equation}
By means of equation (\ref{duhdt}) along with statements (\ref{uptalfab12}) and
(\ref{ghfvhpt}), we have
$$
\frac{\partial \widehat{u}}{\partial t}\in L^{2}({\mathbb R}\times [0,T]),
$$
so that
\begin{equation}
\label{dudtl2}  
\frac{\partial u}{\partial t}\in L^{2}({\mathbb R}\times [0,T]).
\end{equation}
Using the definition of the norm (\ref{142n}) along with 
(\ref{uxtl2}), (\ref{d2udx2l2}) and (\ref{dudtl2}), we derive that
for the function uniquely determined by (\ref{aux}), 
$$
u(x,t)\in W^{1, (4, 2)}({\mathbb R}\times [0, T]).
$$
This means that under the given conditions equation (\ref{aux}) defines a map
$$
\tau_{a, b}: W^{1, (4, 2)}({\mathbb R}\times [0, T])\to
W^{1, (4, 2)}({\mathbb R}\times [0, T]).
$$
Let us establish  that under the stated assumptions  such map is a
strict contraction. We choose arbitrarily
$v_{1, 2}(x,t)\in W^{1, (4, 2)}({\mathbb R}\times [0, T])$. By virtue of the
argument above,
$u_{1, 2}:=\tau_{a, b}v_{1, 2}\in W^{1, (4, 2)}({\mathbb R}\times [0, T])$.
According to formula  (\ref{aux}), we have
$$
\widehat{u_{1}}(p,t)=
$$
\begin{equation}
\label{1aux}
e^{t\{-p^{4}+ibp+a\}}\widehat{u_{0}}(p)+\int_{0}^{t}
e^{(t-s)\{-p^{4}+ibp+a\}}\sqrt{2\pi}\widehat{G}(p)\widehat{f}_{v_{1}}(p,s)ds,
\end{equation}
$$
\widehat{u_{2}}(p,t)=
$$
\begin{equation}
\label{2aux}
e^{t\{-p^{4}+ibp+a\}}\widehat{u_{0}}(p)+\int_{0}^{t}
e^{(t-s)\{-p^{4}+ibp+a\}}\sqrt{2\pi}\widehat{G}(p)\widehat{f}_{v_{2}}(p,s)ds,
\end{equation}
where  $\widehat{f}_{v_{j}}(p,s)$ with $j=1, 2$ denotes the Fourier image of
$F(v_{j}(x,s), x)$ under transform (\ref{ft}). From system
(\ref{1aux}), (\ref{2aux}) we easily obtain that
$$
\widehat{u_{1}}(p,t)-\widehat{u_{2}}(p,t)=
$$
\begin{equation}
\label{u1u2hint0t}  
\int_{0}^{t}
e^{(t-s)\{-p^{4}+ibp+a\}}\sqrt{2\pi}\widehat{G}(p)
[\widehat{f}_{v_{1}}(p,s)-\widehat{f}_{v_{2}}(p,s)]ds.
\end{equation}
Clearly, the upper bound on the norm 
$$
\|\widehat{u_{1}}(p,t)-\widehat{u_{2}}(p,t)\|_{L^{2}({\mathbb R})}\leq
$$
\begin{equation}
\label{int0tetsf12}  
\int_{0}^{t}\|
e^{(t-s)\{-p^{4}+ibp+a\}}\sqrt{2\pi}\widehat{G}(p)
[\widehat{f}_{v_{1}}(p,s)-\widehat{f}_{v_{2}}(p,s)]\|_{L^{2}({\mathbb R})}ds
\end{equation}
holds.
Let us use inequality (\ref{fub}) to derive the upper bound as
$$
\|e^{(t-s)\{-p^{4}+ibp+a\}}\sqrt{2\pi}\widehat{G}(p)
[\widehat{f}_{v_{1}}(p,s)-\widehat{f}_{v_{2}}(p,s)]\|_{L^{2}({\mathbb R})}^{2}=
$$
$$
2\pi \int_{-\infty}^{\infty}e^{-2(t-s)p^{4}}e^{2(t-s)a}|\widehat{G}(p)|^{2}
|\widehat{f}_{v_{1}}(p,s)-\widehat{f}_{v_{2}}(p,s)|^{2}dp\leq 
$$
$$
2\pi e^{2aT}\|\widehat{G}(p)\|_{L^{\infty}({\mathbb R})}^{2}\int_{-\infty}^{\infty}
|\widehat{f}_{v_{1}}(p,s)-\widehat{f}_{v_{2}}(p,s)|^{2}dp\leq 
$$
$$
e^{2aT}\|G(x)\|_{L^{1}({\mathbb R})}^{2}\|F(v_{1}(x,s), x)-F(v_{2}(x,s), x)\|_
{L^{2}({\mathbb R})}^{2}.
$$
Recalling formula (\ref{lip}) gives us
\begin{equation}
\label{lipl2}  
\|F(v_{1}(x,s), x)-F(v_{2}(x,s), x)\|_{L^{2}({\mathbb R})}\leq l
\|v_{1}(x,s)-v_{2}(x,s)\|_{L^{2}({\mathbb R})},
\end{equation}
so that
$$
\|e^{(t-s)\{-p^{4}+ibp+a\}}\sqrt{2\pi}\widehat{G}(p)
[\widehat{f}_{v_{1}}(p,s)-\widehat{f}_{v_{2}}(p,s)]\|_{L^{2}({\mathbb R})}\leq
$$
\begin{equation}
\label{eatgvl2}  
e^{aT}l\|G(x)\|_{L^{1}({\mathbb R})}\|v_{1}(x,s)-v_{2}(x,s)\|_{L^{2}({\mathbb R})}.   
\end{equation}
By means of (\ref{int0tetsf12}) along with (\ref{eatgvl2}),
$$
\|\widehat{u_{1}}(p,t)-\widehat{u_{2}}(p,t)\|_{L^{2}({\mathbb R})}\leq
$$
$$
e^{aT}l\|G(x)\|_{L^{1}({\mathbb R})}\int_{0}^{T}
\|v_{1}(x,s)-v_{2}(x,s)\|_{L^{2}({\mathbb R})}ds.
$$
According to the Schwarz inequality 
\begin{equation}
\label{sch2}  
\int_{0}^{T}\|v_{1}(x,s)-v_{2}(x,s)\|_{L^{2}({\mathbb R})}ds\leq
\sqrt{\int_{0}^{T}\|v_{1}(x,s)-v_{2}(x,s)\|_{L^{2}({\mathbb R})}^{2}ds}\sqrt{T}.
\end{equation}
Hence,
$$  
\|\widehat{u_{1}}(p,t)-\widehat{u_{2}}(p,t)\|_{L^{2}({\mathbb R})}\leq
$$
\begin{equation}
\label{u1hu2hl2}  
e^{aT}l\sqrt{T}\|G(x)\|_{L^{1}({\mathbb R})}\|v_{1}(x,t)-v_{2}(x,t)\|_
{L^{2}({\mathbb R}\times [0, T])}.
\end{equation}
We arrive at
$$
\|u_{1}(x,t)-u_{2}(x,t)\|_{L^{2}({\mathbb R}\times [0, T])}^{2}=\int_{0}^{T}
\|\widehat{u_{1}}(p,t)-\widehat{u_{2}}(p,t)\|_{L^{2}({\mathbb R})}^{2}dt\leq
$$
\begin{equation}
\label{u1u2l2v1v2}
e^{2aT}l^{2}T^{2}\|G(x)\|_{L^{1}({\mathbb R})}^{2}
\|v_{1}(x,t)-v_{2}(x,t)\|_{L^{2}({\mathbb R}\times [0, T])}^{2}.
\end{equation}  
Let us recall (\ref{u1u2hint0t}). Thus,
$$
p^{4}[\widehat{u_{1}}(p,t)-\widehat{u_{2}}(p,t)]=
\int_{0}^{t}
e^{(t-s)\{-p^{4}+ibp+a\}}\sqrt{2\pi}p^{4}\widehat{G}(p)
[\widehat{f}_{v_{1}}(p,s)-\widehat{f}_{v_{2}}(p,s)]ds.
$$
The norm can be trivially estimated from above as
$$
\|p^{4}[\widehat{u_{1}}(p,t)-\widehat{u_{2}}(p,t)]\|_{L^{2}({\mathbb R})}\leq
$$
\begin{equation}
\label{int0tetsf122}  
\int_{0}^{t}\|
e^{(t-s)\{-p^{4}+ibp+a\}}\sqrt{2\pi}p^{4}\widehat{G}(p)
[\widehat{f}_{v_{1}}(p,s)-\widehat{f}_{v_{2}}(p,s)]\|_{L^{2}({\mathbb R})}ds.
\end{equation}
By virtue of inequality (\ref{fub1}) we obtain the estimate from above
$$
\|e^{(t-s)\{-p^{4}+ibp+a\}}\sqrt{2\pi}p^{4}\widehat{G}(p)
[\widehat{f}_{v_{1}}(p,s)-\widehat{f}_{v_{2}}(p,s)]\|_{L^{2}({\mathbb R})}^{2}=
$$
$$
2\pi \int_{-\infty}^{\infty}e^{-2(t-s)p^{4}}e^{2(t-s)a}|p^{4}\widehat{G}(p)|^{2}
|\widehat{f}_{v_{1}}(p,s)-\widehat{f}_{v_{2}}(p,s)|^{2}dp\leq 
$$
$$
2\pi e^{2aT}\|p^{4}\widehat{G}(p)\|_{L^{\infty}({\mathbb R})}^{2}\int_{-\infty}^{\infty}
|\widehat{f}_{v_{1}}(p,s)-\widehat{f}_{v_{2}}(p,s)|^{2}dp\leq 
$$
$$
e^{2aT}\Big\|\frac{d^{4}G}{dx^{4}}\Big\|_{L^{1}({\mathbb R})}^{2}
\|F(v_{1}(x,s), x)-F(v_{2}(x,s), x)\|_{L^{2}({\mathbb R})}^{2}.
$$
Using formula (\ref{lipl2}), we derive
$$
\|e^{(t-s)\{-p^{4}+ibp+a\}}\sqrt{2\pi}p^{4}\widehat{G}(p)
[\widehat{f}_{v_{1}}(p,s)-\widehat{f}_{v_{2}}(p,s)]\|_{L^{2}({\mathbb R})}\leq 
$$
\begin{equation}
\label{etsp2gpeat}
e^{aT}l\Big\|\frac{d^{4}G}{dx^{4}}\Big\|_{L^{1}({\mathbb R})}\|v_{1}(x,s)-v_{2}(x,s)\|_
{L^{2}({\mathbb R})}.
\end{equation}
Let us combine bounds  (\ref{int0tetsf122}), (\ref{etsp2gpeat}) and  (\ref{sch2}). This yields
$$
\|p^{4}[\widehat{u_{1}}(p,t)-\widehat{u_{2}}(p,t)]\|_{L^{2}({\mathbb R})}\leq
$$
\begin{equation}
\label{p2u1hu2hpt}  
e^{aT}\sqrt{T}l\Big\|\frac{d^{4}G}{dx^{4}}\Big\|_{L^{1}({\mathbb R})}
\|v_{1}(x,t)-v_{2}(x,t)\|_{L^{2}({\mathbb R}\times [0, T])}.
\end{equation}
Therefore,
$$
\Big\|\frac{\partial^{4}}{\partial x^{4}}[u_{1}(x,t)-u_{2}(x,t)]\Big\|_
{L^{2}({\mathbb R}\times [0, T])}^{2}=\int_{0}^{T}
\|p^{4}[\widehat{u_{1}}(p,t)-\widehat{u_{2}}(p,t)]\|_{L^{2}({\mathbb R})}^{2}dt\leq     $$
\begin{equation}
\label{u1u2l2v1v22}
e^{2aT}l^{2}T^{2}\Big\|\frac{d^{4}G}{dx^{4}}\Big\|_{L^{1}({\mathbb R})}^{2}
\|v_{1}(x,t)-v_{2}(x,t)\|_{L^{2}({\mathbb R}\times [0, T])}^{2}.
\end{equation}  
By means of (\ref{u1u2hint0t}), we have
$$
\frac{\partial}{\partial t}[\widehat{u_{1}}(p,t)-\widehat{u_{2}}(p,t)]=
$$
$$
\{-p^{4}+ibp+a\}[\widehat{u_{1}}(p,t)-\widehat{u_{2}}(p,t)]+
\sqrt{2\pi}\widehat{G}(p)[\widehat{f}_{v_{1}}(p,t)-\widehat{f}_{v_{2}}(p,t)].
$$
Hence,
$$
\Big\|\frac{\partial}{\partial t}[\widehat{u_{1}}(p,t)-\widehat{u_{2}}(p,t)]
\Big\|_{L^{2}({\mathbb R})}\leq a\|\widehat{u_{1}}(p,t)-\widehat{u_{2}}(p,t)\|_
{L^{2}({\mathbb R})}+
$$
$$      
|b|\|p[\widehat{u_{1}}(p,t)-\widehat{u_{2}}(p,t)]\|_
{L^{2}({\mathbb R})}+
\|p^{4}[\widehat{u_{1}}(p,t)-\widehat{u_{2}}(p,t)]\|_{L^{2}({\mathbb R})}+    
$$
\begin{equation}
\label{ddtu1hptu2hpt}
\sqrt{2\pi}\|\widehat{G}(p)[\widehat{f}_{v_{1}}(p,t)-\widehat{f}_{v_{2}}(p,t)]\|_
{L^{2}({\mathbb R})}.
\end{equation}  
According to  (\ref{u1hu2hl2}), the first term in the right side of
(\ref{ddtu1hptu2hpt}) can be trivially estimated from above by
\begin{equation}
\label{au1hu2hl2}  
age^{aT}\sqrt{T}l\|v_{1}(x,t)-v_{2}(x,t)\|_{L^{2}({\mathbb R}\times [0, T])}
\end{equation}
with $g$ is introduced in (\ref{g}).
We derive the upper bound on the norm as
$$
\|p[\widehat{u_{1}}(p,t)-\widehat{u_{2}}(p,t)]\|_{L^{2}({\mathbb R})}^{2}=
$$
$$
\int_{|p|\leq 1}p^{2}|\widehat{u_{1}}(p,t)-\widehat{u_{2}}(p,t)|^{2}dp+
\int_{|p|>1}p^{2}|\widehat{u_{1}}(p,t)-\widehat{u_{2}}(p,t)|^{2}dp\leq 
$$
$$  
\|\widehat{u_{1}}(p,t)-\widehat{u_{2}}(p,t)\|_{L^{2}({\mathbb R})}^{2}+
\|p^{4}[\widehat{u_{1}}(p,t)-\widehat{u_{2}}(p,t)]\|_{L^{2}({\mathbb R})}^{2}.
$$
Let us use inequalities (\ref{u1hu2hl2}) and  (\ref{p2u1hu2hpt}). They
give the estimate from above on the second term in the right side of
(\ref{ddtu1hptu2hpt}) equal to
\begin{equation}
\label{bu1hu2hl2}  
|b|ge^{aT}\sqrt{T}l\|v_{1}(x,t)-v_{2}(x,t)\|_{L^{2}({\mathbb R}\times [0, T])}.
\end{equation}
By virtue of (\ref{p2u1hu2hpt}), the third term in the
right side of (\ref{ddtu1hptu2hpt}) can be bounded from above by
\begin{equation}
\label{bu1hu2hl3}  
ge^{aT}\sqrt{T}l\|v_{1}(x,t)-v_{2}(x,t)\|_{L^{2}({\mathbb R}\times [0, T])}.
\end{equation}
By means of (\ref{fub}) and (\ref{lipl2}), we obtain that
$$
2\pi\int_{-\infty}^{\infty}|\widehat{G}(p)|^{2}
|\widehat{f}_{v_{1}}(p,t)-\widehat{f}_{v_{2}}(p,t)|^{2}dp\leq
$$
$$
2\pi
\|\widehat{G}(p)\|_{L^{\infty}({\mathbb R})}^{2}\int_{-\infty}^{\infty}
|\widehat{f}_{v_{1}}(p,t)-\widehat{f}_{v_{2}}(p,t)|^{2}dp\leq
$$
$$
\|G(x)\|_{L^{1}({\mathbb R})}^{2}\|F(v_{1}(x,t), x)-F(v_{2}(x,t), x)\|_
{L^{2}({\mathbb R})}^{2}\leq
$$
$$
\|G(x)\|_{L^{1}({\mathbb R})}^{2}l^{2}\|v_{1}(x,t)-v_{2}(x,t)\|_{L^{2}({\mathbb R})}^{2}.  
$$
Thus, the fourth term in the right side of (\ref{ddtu1hptu2hpt}) can be estimated
from above by
\begin{equation}
\label{bu1hu2hl4}  
gl\|v_{1}(x,t)-v_{2}(x,t)\|_{L^{2}({\mathbb R})}.
\end{equation}
Let us combine (\ref{au1hu2hl2}), (\ref{bu1hu2hl2}), (\ref{bu1hu2hl3}) and
(\ref{bu1hu2hl4}). This yields
$$
\Big\|\frac{\partial}{\partial t}[\widehat{u_{1}}(p,t)-\widehat{u_{2}}(p,t)]
\Big\|_{L^{2}({\mathbb R})}\leq
$$
$$
ge^{aT}\sqrt{T}l\{a+|b|+1\}\|v_{1}(x,t)-v_{2}(x,t)\|_
{L^{2}({\mathbb R}\times [0, T])}+
gl\|v_{1}(x,t)-v_{2}(x,t)\|_{L^{2}({\mathbb R})}.
$$
We arrive at the upper bound on the norm as
$$
\Big\|\frac{\partial}{\partial t}(u_{1}(x,t)-u_{2}(x,t))\Big\|_
{L^{2}({\mathbb R}\times [0, T])}^{2}=\int_{0}^{T}
\Big\|\frac{\partial}{\partial t}[\widehat{u_{1}}(p,t)-\widehat{u_{2}}(p,t)]
\Big\|_{L^{2}({\mathbb R})}^{2}dt\leq           
$$
\begin{equation}
\label{ddtu1u2l2}  
2g^{2}l^{2}[e^{2aT}T^{2}\{a+|b|+1\}^{2}+1]
\|v_{1}(x,t)-v_{2}(x,t)\|_{L^{2}({\mathbb R}\times [0, T])}^{2}.
\end{equation}
We recall the definition of the norm (\ref{142n}) and use inequalities
(\ref{u1u2l2v1v2}), (\ref{u1u2l2v1v22}) and (\ref{ddtu1u2l2}). 
A trivial calculation gives us
$$
\|u_{1}-u_{2}\|_{W^{1, (4, 2)}({\mathbb R}\times [0, T])}\leq
$$
\begin{equation}
\label{contr}
gl\sqrt{T^{2}e^{2aT}(1+2[a+|b|+1]^{2})+2}
\|v_{1}-v_{2}\|_{W^{1, (4, 2)}({\mathbb R}\times [0, T])}.
\end{equation}
The constant in the right side of (\ref{contr}) is less than one via
inequality (\ref{qlt}). This means that under the stated assumptions  problem
(\ref{aux}) defines the map
$$
\tau_{a, b}: W^{1, (4, 2)}({\mathbb R}\times [0, T])\to W^{1, (4, 2)}
({\mathbb R}\times [0, T]),
$$
which is a strict contraction. Its unique fixed
point $w(x,t)$ is the only solution of
problem (\ref{h}), (\ref{ic}) in
$W^{1, (4, 2)}({\mathbb R}\times [0, T])$. \hfill\lanbox

\bigskip

\noindent
{\it Proof of Corollary 1.4.} The validity of the statement of our Corollary
follows from the fact that the constant in the right side of bound
(\ref{contr}) is independent of the initial condition (\ref{ic}) (see e.g. ~\cite{EOE20}). Therefore,
problem (\ref{h}), (\ref{ic}) admits a unique solution
$w(x,t)\in W^{1, (4, 2)}({\mathbb R}\times {\mathbb R}^{+})$. We suppose
that $w(x,t)$ is trivial for $x\in {\mathbb R}$ and $t\in {\mathbb R}^{+}$.
This will contradict to the assumption that
$\hbox{supp}\widehat{F(0, x)}\cap \hbox{supp}\widehat{G}$ is a set of
nonzero Lebesgue measure on ${\mathbb R}$. \hfill\lanbox

\bigskip


\centerline{\bf 3. Acknowledgement}

\bigskip

\noindent
V.V. is grateful to Israel Michael Sigal for the partial support by the
NSERC grant NA 7901.


\bigskip

\begin{thebibliography}{99}

\bibitem{ATEYM06}
M. ~Aida, T. ~Tsujikawa, M. ~Efendiev, A. ~Yagi, M. ~Mimura.
{\em Lower estimate of the attractor dimension for a chemotaxis growth system},
J. London Math. Soc. (2), {\bf 74} (2006), no. 2, 453--474.

\bibitem{ABVV10}
N. ~Apreutesei, N. ~Bessonov, V. ~Volpert, V. ~Vougalter.
{\em Spatial structures and generalized travelling waves for an integro-   
differential equation}, Discrete Contin. Dyn. Syst. Ser. B,
{\bf 13} (2010), no. 3, 537--557.

\bibitem{BNPR09}
H. ~Berestycki, G. ~Nadin, B. ~Perthame, L. ~Ryzhik.
{\em The non-local Fisher-KPP equation: travelling waves and steady states},
Nonlinearity, {\bf 22} (2009), no. 12, 2813--2844.

\bibitem{BHN05}
H. ~Berestycki, F. ~Hamel, N. ~Nadirashvili.
{\em The speed of propagation for KPP type problems. I. Periodic 
framework}, J. Eur. Math. Soc. (JEMS), {\bf 7} (2005), no. 2, 173--213. 


\bibitem{BO86}
H. ~Brezis, L. ~Oswald. {\em Remarks on sublinear elliptic equations},
Nonlinear Anal., {\bf 10} (1986), no. 1, 55--64.  

\bibitem{DF06} L. ~Desvillettes,  K. ~Fellner.  {\em Exponential decay 
toward equilibrium via entropy methods for reaction-diffusion equations},
J. Math. Anal. Appl.,  {\bf 319}  (2006), no. 1, 157--176.

\bibitem{DKV20}
S. ~Dipierro, A. ~Karakhanyan, E. ~Valdinoci. {\em A free boundary problem
driven by the biharmonic operator}, Pure Appl. Anal., {\bf 2} (2020), no. 4, 875--942.


\bibitem{DMV05} A. ~Ducrot, M. ~Marion, V. ~Volpert.  {\em Syst\'emes de
r\'eaction-diffusion sans propri\'et\'e de Fredholm},
C. R. Math. Acad. Sci. Paris,  {\bf 340} (2005), no. 9, 659--664.


\bibitem{DMV08} A. ~Ducrot, M. ~Marion, V. ~Volpert.
{\em Reaction-diffusion problems with non-Fredholm operators},
Adv. Differerential Equations, {\bf 13} (2008), no. 11-12, 1151--1192.


\bibitem{E09} M. ~Efendiev.
{\em Fredholm structures, topological invariants and applications}. AIMS Ser. Differ. Equ. Dyn. Syst., {\bf 3}
American Institute of Mathematical Sciences (AIMS), Springfield, MO (2009),  205 pp.

\bibitem{E10} M. ~Efendiev.
{\em Finite and infinite dimensional attractors for evolution equations of mathematical physics}, GAKUTO Internat. Ser. 
Math. Sci. Appl., {\bf 33} $Gakkot{\bar o}sho \ Co.,  \ Ltd., \ Tokyo$ (2010), 239 pp.

\bibitem{E13} M. ~Efendiev.
{\em Evolution equations arising in the modelling of life sciences}, Internat. Ser. Numer. Math., {\bf 163}
Birkhäuser/Springer Basel AG, Basel (2013), 217 pp.

\bibitem{E131} M. ~Efendiev.
{\em Attractors for degenerate parabolic type equations}. Math. Surveys Monogr., {\bf 192} American Mathematical Society, Providence, RI; Real Sociedad Matemática Española, Madrid (2013), 221 pp.

\bibitem{EOE20} M.A. ~Efendiev, M. ~Otani, H.J. ~Eberl.
{\em Mathematical analysis of a PDE-ODE coupled model of mitochondrial 
swelling with degenerate calcium ion diffusion}, SIAM J. Math. Anal., {\bf 52}  
(2020), no. 1, 543--569.


\bibitem{EP07} M. A. ~Efendiev, L.A. ~Peletier.
{\em On the large time behavior of solutions of fourth order parabolic equations and $\epsilon$-entropy of their attractors},
C. R. Math. Acad. Sci. Paris, {\bf 344} (2007), no. 2, 93--96.


\bibitem{EV20} M. ~Efendiev, V. ~Vougalter.
{\em Solvability of some integro-differential equations with drift},
Osaka J. Math., {\bf 57} (2020), no. 2, 247--265.


\bibitem{EV22} M. ~Efendiev, V. ~Vougalter.
{\em Solvability of some integro-differential equations with drift and
superdiffusion}, J. Dynam. Differential Equations, {\bf 36} (2024), no. 1,
353--373.

\bibitem{EV25} M. ~Efendiev, V. ~Vougalter.
{\em On the well-posedness of some model arising in the mathematical biology},
Math. Methods Appl. Sci., {\bf 48} (2025), no. 3, 3670--3681.

\bibitem{GVA06} S. ~Genieys, V. ~Volpert, P. ~Auger.
{\em Pattern and waves for a model in population dynamics with nonlocal 
consumption of resources}, Math. Model. Nat. Phenom., {\bf 1} (2006), no. 1, 
65--82.

\bibitem{GHLP22} V. ~Giunta, T. ~Hillen, M. A. ~Lewis, J.R. ~Potts.
{\em Local and global existence for nonlocal multispecies advection-diffusion models}, 
SIAM J. Appl. Dyn. Syst., {\bf 21} (2022), no. 3, 1686--1708.

\bibitem{GK18} J.R.G. ~Granada, V.A. ~Kovtunenko.
{\em Entropy method for generalized Poisson-Nernst-Planck equations}, Anal.
Math. Phys., {\bf 8} (2018), no. 4, 603--619.  


\bibitem{K64} M.A. ~Krasnosel'skii.
{\em Topological methods in the theory of nonlinear integral equations}.
The Macmillan Company, New York (1964), XI, 395 pp.


\bibitem{LL97} E.H. ~Lieb, M.~Loss.
{\em Analysis. Grad. Stud. Math.}, {\bf 14}, American
 Mathematical Society, Providence, RI (1997), XVIII, 278 pp.

\bibitem{PT01} L.A. ~Peletier, W.C. ~Troy.
{\em Spatial patterns}. Higher order models in physics and mechanics.
Progr. Nonlinear Differential Equations Appl., {\bf 45}.
Birkh\"auser Boston, Inc., Boston, MA (2001), 341 pp.


\bibitem{VV130} V. ~Volpert, V. ~Vougalter. {\em Emergence and propagation
of patterns in nonlocal reaction-diffusion equations arising in the theory
of speciation.} 
Lecture Notes in Math., {\bf 2071} (2013), Springer, Heidelberg, 331--353.


\bibitem{VV210} V. ~Vougalter, V. ~Volpert. {\em Existence of solutions for some 
non-Fredholm integro-differential equations with the bi-Laplacian}, Math. Methods Appl. Sci., 
{\bf 44} (2021), no. 1, 220--231.



\bibitem{VV21} V. ~Vougalter, V. ~Volpert. {\em Solvability of some integro-
differential equations with anomalous diffusion and transport},
Anal. Math. Phys., {\bf 11} (2021), no. 3, Paper No. 135, 26 pp.     

  
\end{thebibliography}
\end{document}